\documentclass[12pt, A4paper]{article}

\usepackage{graphicx} 
\usepackage[utf8]{inputenc}
\usepackage{amsthm,amsmath,amssymb,enumerate,comment}
\usepackage{xcolor,colortbl}
\usepackage{hyperref}

\newtheorem{conjecture}{Conjecture}
\newtheorem{theorem}{Theorem}[section]
\newtheorem{lemma}[theorem]{Lemma}
\newtheorem{proposition}[theorem]{Proposition}

\theoremstyle{definition}

\newtheorem{remark}{Remark}

\newcommand{\ZZ}{\mathbb{Z}}

\newcommand{\arxiv}[2]{\href{https://arxiv.org/abs/#1}{\texttt{arXiv:#1}} \texttt{[#2]}}

\title{Hamilton cycles in bicirculants\\ with non-coprime connection sets}

\usepackage{authblk}

\author[1]{Simona\,Bonvicini}
\author[2]{Tomaž Pisanski}
\author[3,4]{Arjana Žitnik}
\affil[1]{Dipartimento di Scienze Fisiche, Informatiche e Matematiche, \mbox{Universit\`a di Modena e Reggio Emilia, via Campi 213/b, 41126 Modena, Italy}}
\affil[2]{University of Primorska, IAM and FAMNIT, Koper, Slovenia}
\affil[3]{University of Ljubljana, Faculty of Mathematics and Physics, Ljubljana, Slovenia}
\affil[4]{Institute of Mathematics, Physics and Mechanics, Ljubljana, Slovenia}

\date{\today}

\begin{document}

\maketitle

\begin{abstract}
A \emph{bicirculant} is a regular graph that admits an automorphism having two vertex-orbits of the same size. 
A  bicirculant can be described as follows. Given an integer $m \ge 1$ and sets $R, S, T \subseteq \ZZ_m$ such that $R=-R$, $T=-T$, $0 \not\in R \cup T$ and $0 \in S$, the graph $B(m;R,S,T)$ has vertex set $V=\{u_0,\dots,u_{m-1},v_0,\dots,v_{m-1}\}$ and edge set 
$E=\{u_iu_{i+j}| \ i \in\ZZ_m, j \in R\} \cup 
      \{v_iv_{i+j}| \ i \in\ZZ_m, j \in T\} \cup
      \{u_iv_{i+j}| \ i \in\ZZ_m, j \in S\}.$
Some of the well known bicirculants include the generalized Petersen graphs and Cayley graphs of dihedral groups. 

In 1983, Brian Alspach proved that the only non-hamiltonian generalized Petersen graphs are $G(m,2)$ where $m \equiv 5 \pmod 6$.  In a recent paper we conjectured  that this is the only exception among regular connected bicirculants of degree at least two, and we have verified the conjecture for bicirculants with $|S|\le 3$. Results of Alspach, Chen and Dean imply that the conjecture also holds for all bicirculants with $m$ even and $\gcd(m,S)=1$.
In this paper, we show that the conjecture holds for all bicirculants with $|S| \geq 4$ and $\gcd(m,S)>1$ such that
\begin{itemize}
\item $R\cup T$ contains an element that is not coprime to $\gcd(m,S)$ and, if $m$ is odd,  is also not divisibile by $\gcd(m, S)$ or 
\item $R \cap T$ contains an element that is coprime to $\gcd(m,S)$.
\end{itemize}
\end{abstract}


\textbf{Math Subj Class (2020)}: 
05C45, 
05C25, 
05C76,  
05C70,  
05E18. 


\section{Introduction}

A \emph{bicirculant} is a regular graph  that admits an automorphism having two vertex-orbits of the same size \cite{Pi2007}.  
A typical bicirculant can be described as follows; see for example \cite{MALNIC2007891}. Given an integer $m \ge 1$ and sets $R,S,T \subseteq \ZZ_m$
such that $R=-R$, $T=-T$,  
$0 \not\in R \cup T$   and $0 \in S$, 
the graph $B(m;R,S,T)$ has vertex set 
$V=V_1 \cup V_2$, where $V_1=\{u_0,\dots,u_{m-1}\}$ and $V_2=\{v_0,\dots,v_{m-1}\}$, and edge set 
$$E=\{u_iu_{i+j}| \ i \in\ZZ_m, j \in R\} \cup 
      \{v_iv_{i+j}| \ i \in\ZZ_m, j \in T\} \cup
      \{u_iv_{i+j}| \ i \in\ZZ_m, j \in S\}.$$ 
It is easy to see that the mapping $\alpha:V \to V$, defined by $\alpha(u_i)=u_{i+1}$, $\alpha(v_i)=v_{i+1}$
is an automorphism of $B(m;R,S,T)$, having two vertex-orbits of the same size.
We consider only regular bicirculants, so $|R|=|T|$.

The vertices belonging to $V_1$ and $V_2$ will be called \emph{outer vertices} and  \emph{inner vertices}, respectively.
The subgraph induced by the  outer (inner) vertices is a circulant that is described by the set $R$  
($T$,  respectively).  It is called the \emph{outer rim} (the \emph{inner rim}, respectively).
We distinguish three types of edges:  outer edges, inner edges and  spokes.
The \emph{outer edges} are edges incident to two outer vertices; 
outer edges $u_iu_{i+a}$, $i \in \ZZ_m$, $a \in R$, are called   \emph{outer edges of type} $a$.
The \emph{inner edges} are edges incident to two inner vertices;   inner edges $v_iv_{i+b}$, $i \in \ZZ_m$, $b \in T$, are called \emph{inner edges of type} $b$. The \emph{spokes} are edges that connect an outer vertex to an inner vertex;  spokes $u_iv_{i+c}$, $i \in \ZZ_m$, $c \in S$, are called  \emph{spokes of type} $c$.  
The spokes induce a bipartite subgraph of $B(m;R,S,T)$, denoted by $H(m;S)$. Observe that $H(m;S)$ is a bicirculant of the form $B(m;\emptyset,S,\emptyset)$. Such bicirculants are called \emph{cyclic Haar graphs} \cite{HlMaPi2002}.

Bicirculants include many  well known families of graphs such as the generalized Petersen graphs and more generally
$I$-graphs $I(m;j,k)$, which are  bicirculants $B(m;R,S,T)$ with $m \ge 3$, $R=\{j,-j\}$, $T=\{k,-k\}$ and $S=\{0\}$, where $0 < j,k < m/2$.  
The generalized rose window graphs  \cite{wilsonRW, BoPiZi2025}  are $4$-valent bicirculants that can be obtained by adding an additional set of spokes to an $I$-graph. The generalized Taba\v cjn graphs \cite{Taba1,Taba2} are pentavalent bicirculants whose spokes induce   cubic cyclic Haar graphs. Also Cayley graphs of dihedral groups are bicirculant graphs  \cite{MALNIC2007891}.

A natural question concerning bicirculants is to determine which of them admit a hamilton cycle. This question was first studied by Brian Alspach for generalized Petersen graphs, and he gave a complete classification of the hamiltonian generalized Petersen graphs \cite{Al11983}. Hamiltonicity has since been studied for several other families of bicirculants, including $I$-graphs \cite{BoPi2017} and Cayley graphs of dihedral groups \cite{AlChDe2010,AlZh1989,Pastine2012}. In our previous work \cite{BoPiZi2025}, we showed that every generalized rose window graph is hamiltonian. In all these families, the only non-hamiltonian examples are those already identified by Alspach. This led us to formulate the following conjecture \cite{BoPiZi2025}.

\begin{conjecture} \label{conject}
Every connected bicirculant, except for $K_2$ and the generalized Petersen graphs $G(m, 2)$ with $m\equiv 5\pmod 6$, is hamiltonian. 
\end{conjecture}

In our subsequent papers \cite{BoPiZi2025_2ndpaper,BoPiZi2025_3rdpaper}, we obtained several partial results towards the conjecture, including that all bicirculants $B(m;R,S,T)$ with $|S| \le 3$ are hamiltonian, except for the graphs found by Alspach.
Alspach, Chen and Dean \cite{AlChDe2010} proved that every connected cyclic Haar graph $H(m;S)$ is hamiltonian for $m$ even and has a hamilton path for $m$ odd.
It follows that for $m$ even and $\gcd(m,S)=1$, every bicirculant  $B(m;R,S,T)$ is hamiltonian, since it contains as a subgraph the connected cyclic Haar graph $H(m;S)$. However, the question whether cyclic Haar graphs $H(m;S)$ are hamiltonion for $m$ odd is still wide open for $|S| \ge 4$. In fact, if all cyclic Haar graphs $H(m;S)$ are hamiltonian also for $m$ odd, $m>1$, then all bicirculants (except for the known exceptions) are hamiltonian by \cite[Theorem 1.3]{BoPiZi2025_3rdpaper}. In this paper, we consider the remaining cases with $|S| \ge 4$ and $\gcd(m,S) >1$. The following are our main results.

\begin{theorem} \label{th_main}
    Let $G=B(m, R,S,T)$ be a connected bicirculant graph with $|S|\ge 4$ and $\gcd(m, S)>1$. Assume that $R\cup T$ contains an element that is not coprime to $\gcd(m, S)$ and, if $m$ is odd, is not divisible by $\gcd(m, S)$. Then $G$ is hamiltonian. 
\end{theorem}

\begin{theorem}\label{th_a=b}
Let $G=B(m; R,  S,  T)$ be a connected bicirculant graph with  $|S|\ge 4$ and $\gcd(m,S)>1$. Assume that $R \cap T$ contains an element that is coprime to $\gcd(m,S)$. Then $G$ is hamiltonian.
\end{theorem}

Our results apply to bicirculant graphs $B(m;R,S,T)$ for which $\gcd(m,S)>1$. In this case the subgraph $H(m;S)$ is disconnected. We introduce constructions that produce a hamilton cycle in the whole graph by appropriately joining selected paths from the connected components of  $H(m;S)$. In contrast to our previous work \cite{BoPiZi2025_2ndpaper,BoPiZi2025_3rdpaper}, where we required the connected components of $H(m;S)$ to be hamiltonian, the present approach does not rely on this assumption. This allows us to obtain results beyond those covered in our previous papers on this topic. Only a small proportion of bicirculant graphs remain for which it is not known whether they are hamiltonian; these cases are discussed in the final section.


\section{Preliminary results}\label{sec_preliminary}

In this section, we describe the state of the art on the existence of hamilton cycles in bicirculants by reviewing the results that are involved in the proofs of the main theorems of this paper. We also derive an auxiliary lemma describing particular hamilton paths  
in bicirculants that will be useful later. At the end of the section, we review some basic properties of bicirculants.

\smallskip

One of the first results on the existence of hamilton cycles in bicirculants concerns Cayley graphs of valence $3$ on dihedral groups: Alspach and Zhang proved that they are hamiltonian  \cite[Theorem 3.1]{AlZh1989}.  A stronger property on the hamilton-connectivity, or laceability,  for Cayley graphs on generalized dihedral groups was proved later by Alspach, Chen and Dean in  \cite{AlChDe2010}.  We recall that a bicirculant graph $B(m;R,S,T)$ with $R=T$ is a Cayley graph  on the dihedral group of order $2m$. In particular, cyclic Haar graphs are Cayley graphs on dihedral groups of form $B(m; \emptyset, S,\emptyset)$.
 The following result thus follows directly from  
\cite[Theorem 1.8]{AlChDe2010}. 

\begin{theorem}
\label{th_generalized_dihedral}
Let $G=B(m;R,S,R)$ be a connected bicirculant graph with $m$ even and valency at least $3$. Then $G$ is hamilton-connected, unless it is bipartite, in which case it is hamilton-laceable. In particular, $G$ is hamiltonian.
\end{theorem}

It follows that all connected cyclic Haar graphs $H(m;S)$ with $m$ even are hamiltonian. Unfortunately, Theorem \ref{th_generalized_dihedral} does not cover cyclic Haar graphs $H(m;S)$ for  $m$ odd.
However, it follows from \cite[Lemma 2.9]{AlChDe2010} that such graphs have hamilton paths. We restate the lemma in the language of bicirculants.

\begin{lemma}
\label{lem_ham_path_cyclic_Haar}
Every connected cyclic Haar graph has a hamilton path.
\end{lemma}

The work of Pastine and Jaume \cite{Pastine2012} 
provides additional cases, also for $m$ odd. In particular, it follows from their work that every bicirculant $B(m;R,S,R)$ with $R\neq\emptyset$ is hamiltonian.

We next state results from our previous work on hamilton cycles in bicirculant graphs \cite{BoPiZi2025,BoPiZi2025_2ndpaper,BoPiZi2025_3rdpaper} needed for this paper.
In  \cite{BoPiZi2025_2ndpaper} we proved that several large classes of bicirculant graphs are hamiltonian. In particular, we proved the following results.

\begin{lemma}[\protect{\cite[Lemma 3.1]{BoPiZi2025_2ndpaper}}]\label{lemma_bicirculant_m<5}
All connected bicirculant graphs of order $2m$ with $m\le 5$ are hamiltonian, except for the  complete graph $K_2$ and the 
Petersen graph $G(5,2)$. 
\end{lemma}

\begin{theorem}[\protect{\cite[Theorem 1.3]{BoPiZi2025_2ndpaper}}]\label{th_m/2_2ndpaper}
Let  $G=B(m; R,  S, T)$ be a connected bicirculant graph with $m$ even,  $m\ge 4$, $|S|\ge 3$ and $m/2\in R, T$. Then  $G$ is hamiltonian.
\end{theorem}

\begin{theorem}[\protect{\cite[Theorem 1.4]{BoPiZi2025_2ndpaper}}]\label{th_main_2ndpaper}
Let $G=B(m; R, S, T)$ be a connected bicirculant graph with $m\ge 3$, $|S|\ge 3$.
Assume that the connected components of $H(m; S)$ are hamiltonian and that $R\cup T$ contains an element that is not coprime to $\gcd(m, S)$ in the case where $\gcd(m, S)>1$. Then $G$ is hamiltonian.
\end{theorem}

In  \cite{BoPiZi2025_3rdpaper},  we extended the results of \cite{BoPiZi2025_2ndpaper} by constructing   hamilton cycles also for other families of bicirculants,  which are described in the  following theorems.

\begin{theorem}[\protect{\cite[Theorem 1.1]{BoPiZi2025_3rdpaper}}]\label{th_bicirculant_s3_3rdpaper}
A connected bicirculant graph whose vertices are incident to $s\le 3$ spokes is hamiltonian, except for the complete graph
$K_2$ and the generalized Petersen graphs  $G(m,2)$ where $m \equiv 5 \pmod 6$. 
\end{theorem}

\begin{theorem}[\protect{\cite[Theorem 1.2]{BoPiZi2025_3rdpaper}}]\label{th_bicirculant_m/gcd_3rdpaper}
Let   $G=B(m; R, S, T)$ be a connected bicirculant graph with $|S|\ge 4$ and $m/\gcd(m, S)$ even. Then $G$  is hamiltonian.
\end{theorem}

For the constructions in the next section, we will also need the following lemma. We will use the notation
$xP y$ to denote the subpath of a path $P$, starting at vertex $x$ and ending at vertex $y$.

\begin{lemma}\label{lemma_path_1}
Let $H=H(m; S)$ be a connected cyclic Haar graph with 
$|S|\ge 4$.  Then the graph $H$ is hamiltonian, or there exists a hamilton path $P$ in $H$ connecting an inner vertex $v_s$ to the vertex $u_0$,  with the property that at least one of the neighbors of $v_s$,  say $u_z$, is such that the subpaths  $v_s\,P\,u_z$, $u_z\,P\,u_0$ contain at least $6$ and $7$ vertices, respectively.
\end{lemma}

\begin{proof}
If $m\le 5$, then $H$ is hamiltonian by Lemma \ref{lemma_bicirculant_m<5}. We now assume that $m >5$.
Denote the degree of $H$ by $d$, so $d=|S| \ge 4$.
We know that every cyclic Haar graph has a hamilton path by Lemma \ref{lem_ham_path_cyclic_Haar}. 
Moreover, every cyclic Haar graph is vertex-transitive and bipartite. Therefore we can assume that $H$ has a hamilton path $P$ with terminus $u_0$,  and the origin is an inner vertex, say $v_s$.
If $u_0$ and $v_s$ are adjacent, then $H$ is hamiltonian and the assertion follows. 

Let us consider the case where $u_0$ and $v_s$ are not adjacent. Since $d\ge 4$, the vertex $v_s$ has at least $4$ distinct neighbours that are internal vertices of $P$.  
We set $P=(v_s, u_{x_1}, v_{y_1},\ldots, u_{x_i}, v_{y_i},\ldots, u_{x_t}, v_{y_t}, u_0)$, with $t\ge 5$ since $m>5$.  
If at least one of the vertices $u_{x_2}, u_{x_t}, u_{x_{t-1}}$
is not adjacent to $v_s$, then we can find a vertex $u_z$ satisfying the required property. 
Assume that the vertex $v_s$ is adjacent to  
each of $u_{x_2}, u_{x_t}, u_{x_{t-1}}$.
If $d>4$, then we can find at least one vertex $u_z$ with the required property.

We now consider the worst case, that is, $d=4$ and $v_s$ is adjacent to 
each of $u_{x_2}, u_{x_t}, u_{x_{t-1}}$.
Let us consider the neighbours of $v_{y_{t-1}}$.
If $v_{y_{t-1}}$ is adjacent to $u_0$, then $H$ is hamiltonian: removing the edge $u_{x_t}v_{y_{t-1}}$ from $P$ and adding the edges $u_{x_t}v_s$, $u_0v_{y_{t-1}}$ produces a hamilton cycle in $H$.
If $v_{y_{t-1}}$ is adjacent to $u_{x_j}$ 
for some $j \in \{2,\dots, t-3\}$, 
then we remove the edge $u_{x_t}v_{y_{t-1}}$ from $P$, add the edge $u_{x_t}v_s$, and find a hamilton path from $v_{y_{t-1}}$ to $u_0$ with the vertex  $u_{x_j}$ playing the role of the vertex $u_z$, since $1<j<t-2$.

If none of the vertices $u_0$ and $u_{x_j}$ with $1<j<t-2$ is adjacent to $v_{y_{t-1}}$, it means that $v_{y_{t-1}}$ is adjacent to $u_{x_1}$ and $u_{x_{t-2}}$.
In this case, we remove the edge $u_{x_1}v_s$ and the edges in the subpath $(u_{x_t},v_{y_{t-1}}, u_{x_{t-1}})$ from $P$; we add the path $(u_{x_t}, v_s, u_{x_{t-1}})$ and the edge $v_{y_{t-1}}u_{x_1}$, and find a hamilton path from $v_{y_{t-1}}$ to $u_0$ with the vertex $u_{x_{t-2}}$ playing the role of the vertex $u_z$.
\end{proof}

We end this section by reviewing some properties  that will be used in the paper, see \cite{BoPiZi2025_2ndpaper} and \cite{GT}. 
We recall that we always assume that $0 \in S$ in a bicirculant  $B(m;R,S,T)$.

\begin{proposition}  \label{prop:connected}
A bicirculant graph $B(m;R,S,T)$  is connected if and only if $\gcd(m,R,S,T)=1$. 
\end{proposition}

In the case that a bicirculant graph is disconnected, it is composed of isomorphic connected components. We will use the following notation. Let $A=\{a_1,a_2,\dots,a_\ell\}$ be a set and let $i$ be an integer. We define $A-i=\{a-i \ | \ a \in A \}$ and $A/i= \{a/i \ | \ a \in A \}$. Moreover,  the notation $\gcd(A)$ should be understood as $\gcd(a_1,a_2,\dots,a_\ell)$ and $\gcd(i,A)$ as $\gcd(i,\gcd(A))$.

\begin{proposition}\label{pro:d_connected_components} 
Let $G = B(m;R,S,T)$. Suppose $\delta = \gcd(m,R,S,T)>1$. 
Then $G$ is a disjoint union of $\delta$   isomorphic graphs $G_0, \dots G_{\delta-1}$ such that $u_i \in G_i$, $i = 0,\dots,\delta-1$. Moreover, each graph $G_i$  is connected and isomorphic to the graph $B(m/\delta;R/\delta,S/\delta,T/\delta)$. 
\end{proposition}

Different parameters may give the same bicirculants. In particular, it is easy to verify that the following result holds.

\begin{lemma}  \label{lemma:iso}
Let $G = B(m;R,S,T)$ be a bicirculant and let $c \in S$. Then the graph $B(m;R,S-c,T)$ is isomorphic to the graph $G$.
\end{lemma}

In the rest of the paper, we will denote a bicirculant $B(m;\{a,-a\},S,\{b,-b\})$ simply by $B(m;a,S,b)$.

\section{From hamilton paths in $H(m; S)$ to a hamilton cycle in $B(m; a, S, b)$}
\label{sec_construction}

In this section, we consider connected bicirculant graphs $B(m; a, S, b)$ with  $|S|\ge 4$, at least one of the integers $a, b$ not coprime to $\gcd(m,S)$ and $a, b\ne m/2$.  Note that in this case $m\ge 4$ since $|S|\ge 4$. 
We present different constructions to obtain a hamilton cycle in a graph $B(m; a, S, b)$, starting from  hamilton paths in the connected components of $H(m; S)$ whose existence is guaranteed by  Lemma \ref{lem_ham_path_cyclic_Haar}. 
For our purposes, we will apply the constructions to bicirculant graphs with $|S|\ge 4$, but they can also be applied to graphs with $|S|\ge 2$ provided that the number of vertices in the components of $H(m; S)$ is not too small. 

Our constructions rely on the fact that the graphs considered in this section admit a uniform representation.
Uniform representations have been introduced in \cite{BoPiZi2025_2ndpaper} and formally defined in \cite{BoPiZi2025_3rdpaper}. 
As remarked in \cite[Section 3]{BoPiZi2025_3rdpaper}, we can find a uniform representation for a graph $B(m;a,S,b)$ when at least one of the integers $a, b$ is not coprime to $\gcd(m, S)$. More specifically, let $x\in\{\pm a, \pm b\}$ such that $\gcd(m, S, x)>1$. We set $\lambda=\lambda(x)=\gcd(m, s, x)-1$, $\mu=\mu(x)=\gcd(m, S)/\gcd(m, S, x)-1$, and say that $B(m;a,S,b)$ has a $(\lambda, \mu)$-\emph{uniform representation}. A uniform representation of $B(m;a,S,b)$ allows us to describe the structure of the graph as follows.
\smallskip

Let $G=B(m; a, S, b)$ be a connected bicirculant graph with  $|S|\geq 4$, $\gcd(m, S)>1$,   $a,b\ne m/2$, and at least one of the integers $a$, $b$, say $b$,   not coprime to $\gcd(m,S)$.  We denote with $H=H(m; S)$  the cyclic Haar graph induced by the spokes of $G$, and with  $K$ the subgraph obtained from $G$ by removing the edges of type $a$. The connected components of $K$ form a partition of the components of $H$: each component of $K$ consists of $\gcd(m, S)/\gcd(m, S, b)$ components of $H$, which are joined by edges of type $b$.  We set $\lambda=\lambda(b)=\gcd(m, S, b)-1$, $\mu=\mu(b)=(\gcd(m, S)/\gcd(m, S, b))-1$ and denote with $K_0,\ldots, K_{\lambda}$ the connected components of $K$; by the assumptions we have $\lambda >0$.

For $0\le j\le\lambda$, we denote with $H_{i,j}$, with $0\le i\le\mu$, the connected components of $H$ that are contained in $K_j$.
Moreover, we define $H_{0,0}$ to be the component of $H$ containing the vertex $u_0$; consequently, $K_0$ is the component of $K$ containing $u_0$.

The graph $G$ can be represented as follows: arrange the components $H_{i,j}$ with the same index $i$ in the same row $i$ and the components $H_{i,j}$ with the same index $j$ in the same column $j$. In this setting, the outer vertices in $H_{i, j}$ are adjacent to the outer vertices in $H_{i, j+1}$ for $0\le j\le\lambda-1$, and the inner vertices in $H_{i, j}$ are adjacent to the inner vertices in $H_{i+1, j}$ for $0\le i\le\mu$ (the subscript $i$ is considered modulo $\mu+1$).
We also have that the outer vertices of $K_0$ are adjacent to the outer vertices of $K_{\lambda}$.
Due to this arrangement, the component $K_j$ is called  `the vertical $j$-th component' of $G$, and the representation of $G$ is called a $(\lambda, \mu)$-\emph{uniform representation} of $G$, or simply a uniform representation of $G$, because each vertical component $K_j$ consists of $(\mu+1)$ components of $H(m; S)$. 
In \cite{BoPiZi2025_2ndpaper} we proved the following result for bicirculant graphs having a uniform representation. 

\begin{lemma}[\protect{\cite[Lemma 3.2]{BoPiZi2025_2ndpaper}}]\label{lem_BoPiZi2025_2ndpaper}
Let $G=B(m; a, S, b)$ be a connected bicirculant graph with  $m\ge 2$, $|S|\geq 2$, $\gcd(m, S)>1$,  $a,b\ne m/2$,  and $\lambda=\gcd(m, S, b)-1 >0$.  
Assume that the connected components of $H(m; S)$ are hamiltonian.  Then $G$ is hamiltonian.
\end{lemma}

The constructions that we describe in this section improve upon the result of Lemma \ref{lem_BoPiZi2025_2ndpaper} by dropping the assumption on the existence of a hamilton cycle in the components of $H(m; S)$. Instead, we construct a hamilton cycle in $G$ from hamilton paths in the components of $H(m; S)$.
The constructions are described in Lemmas \ref{pro_method1_case}, \ref{pro_method1_case2} and \ref{pro_method1_case3}, which depend on the parity of the parameters $\lambda$ and $\mu$. 
\smallskip

When the order of the connected components of $H(m;S)$ is small, i.e. $m/\gcd(m,S)$ $\le 5$, and $|S| \ge 2$, the components of $H(m;S)$ are  hamiltonian by Lemma \ref{lemma_bicirculant_m<5}. In this case, Lemma \ref{lem_BoPiZi2025_2ndpaper} can be applied to obtain a hamilton cycle in the whole graph. 
If the components of $H(m,S)$ are not hamiltonian, then we may assume that their order is more than $2 \cdot 5$ and that they contain a hamilton path as described in  Lemma \ref{lemma_path_1}. 
We will use the  notation and terminology from \cite{BoPiZi2025_2ndpaper, BoPiZi2025_3rdpaper}
with appropriate modifications due to the case being studied. In detail, let $P$ be a hamilton path in $H_{0,0}$ 
with origin $v_s$, terminus $u_0$, and suppose that there exists a neighbor $u_z$ of $v_s$ such that the subpaths $v_s\,P\,u_z$ and $u_z\,P\, u_0$ contain at least $6$ and $7$ vertices, respectively.  
We  denote with $v_h$, $v_k$ the vertices adjacent to $u_z$ in $P$. For convenience of notation, the remaining vertices of $P$ will be simply denoted with $u_i$, $v_i$ instead of $u_{x_i}$, $v_{y_i}$. More specifically, we set 

$$P=(v_s, u_t, v_t, u_{t-1},  \ldots,   v_h,  u_z,  v_k,  u_{k},\ldots,  
v_1, u_1,  v_0 , u_0).$$

\noindent 
Since we have assumed that the paths $v_s\,P\,u_z$ and $u_z\,P\,u_0$ of $P$ contain at least $6$ and $7$ vertices, we can assume that the vertices $v_s$, $u_t$, $v_t$, $u_{t-1}$, $v_h$, $u_z$ are pairwise distinct, while $k\ge 2$ in the subpath $u_z\,P\,u_0$.

The connected components of $K$ -- and also of $H$ -- contain copies of $P$ which, for our purposes, we will describe as follows.
For $0\leq i\leq\mu$ and $0\leq j\leq\lambda$, we denote by $P_{i,j}$ the path obtained from $P$ by adding $ib+ja\pmod m$ to the subscripts of the vertices in $P$; $P_{i,j}$ is the copy of the path $P$ in the component $H_{i,j}$ of $H$, and $P_{0,0}$ corresponds to the path $P$. 
Accordingly, the vertices in $P_{i,j}$ will be denoted with $u^{i,j}_x$, $v^{i,j}_x$, where $u^{i,j}_x=u_{x+ib+ja}$, $v^{i,j}_x=v_{x+ib+ja}$.  The outer vertices $u^{i,j}_x$ in $P_{i,j}$ are adjacent to the outer vertices $u^{i,j+1}_x$ in $P_{i, j+1}$ by the edges of type $a$  for $0\le j\le\lambda-1$; the inner vertices
$v^{i,j}_x$ in $P_{i,j}$ are adjacent to the inner vertices $v^{i+1,j}_x$ in $P_{i+1, j}$ by the edges of type $b$ for $0\le i\le\mu-1$. 
For simplicity of notation, we will denote with $u_x\,P_{i,j}\,v_y$, or with  $u_x\,P\,v_y$ in $H_{i,j}$, the path corresponding to the subpath $u_x\,P\,v_y$ from the vertex $u_x$ to the vertex $v_y$ of $P$; 
it is understood that the vertices $u_x$, $v_y$ in the notation $u_x\,P_{i,j}\,v_y$, or in the notation $u_x\,P\,v_y$ in $H_{i,j}$,  correspond to the vertices $u^{i,j}_x$, $v^{i,j}_y$.

\smallskip

The proofs of Lemmas \ref{pro_method1_case}--\ref{pro_method1_case3} 
consist of selecting a collection of special paths in each component of $H(m;S)$ and joining them to form a hamilton cycle in $G$. The selected paths and the way they are connected can be readily seen from the figures that accompany the proofs; see Figures \ref{fig_ham_cycle_lambda_mu_odd}--\ref{fig_ham_cycle_even}. 
In the figures,  the components $H_{i,j}$ of $H(m; S)$ are  
represented as (underlying) circles that are arranged in a rectangular grid consisting of $\lambda+1$ vertical columns with $\mu+1$ components of $H(m; S)$; vertical columns are the vertical components $K_j$, with $0\le j\le\lambda$, which define a uniform representation of $G$;
the path $P$ and its copies lie in the (underlying) circles. The curves that appear in the figure at each $H_{i,j}$ represent the copies of the subpaths of $P$ that are used in the proof of Lemmas  \ref{pro_method1_case}--\ref{pro_method1_case3}. In particular,
straight line segments represent edges, while curved segments represent paths that might have more than one edge.

\smallskip

We now introduce the first construction of this section.

\begin{lemma}\label{pro_method1_case}
Let $G=B(m; a, S, b)$ be a connected bicirculant graph with 
$|S|\geq 4$ and 
$a,b\ne m/2$ 
that has a $(\lambda, \mu)$-uniform representation such that $\lambda=\lambda(b)=\gcd(m, S, b)-1 >0$ and
$\mu=\mu(b)=\gcd(m, S)/\gcd(m, S, b)-1\ge 0$.
Assume that $\lambda$ and $\mu$ are both odd. Then $G$ is hamiltonian.
\end{lemma}

\begin{proof}
If the connected components of $H(m;S)$ are hamiltonian, then the assertion follows from Lemma \ref{lem_BoPiZi2025_2ndpaper}. Now assume that this is not the case. 
The uniform representation of $G$ allows us to describe it as detailed above and, in what follows, we will use this representation to construct a hamilton cycle in $G$.
By Lemma \ref{lemma_path_1} we can assume that $H_{0,0}$ contains a hamilton path $P$  described earlier in this section.

We give a list of  subpaths of $P$ that will be
employed in the construction of a hamilton cycle in $G$. 
We give the list for $\lambda \ge 1$ and $\mu> 1$, and also for the case where $\lambda=\mu=1$.
The case with $\lambda>1$ and $\mu=1$ can be obtained from the case with $\lambda=1$ and $\mu>1$ by exchanging the inner vertices with the outer vertices. More precisely, by Lemma \ref{lemma_path_1} we have a path $P$ from $u_s$ to $v_0$ with the property that at least one of the neighbors of $u_s$, say $v_z$, is such that the subpaths $u_sPv_z$, $u_zPu_0$ contain at least $6$ and $7$ vertices, respectively. In the construction, we use the same indices as in the case  $\lambda \ge 1$ and $\mu> 1$, but with inner and outer vertices exchanged. We also swap $\lambda$ and $\mu$ and $H_{i,j}$ with $H_{j, i}$ everywhere.

For the construction of a hamilton cycle in $G$ when  $\lambda=\mu=1$, we consider a copy of the path $P$ from $v_s$ to $u_0$ in each component $H_{i, j}$ of $H$,
and connect the path $P_{i,0}$ to the path $P_{i,1}$ by adding the edge $u^{i,0}_0u^{i,1}_0$, with $i=0,1$, and also connect the path $P_{0,j}$ to the path $P_{1,j}$ through the edge $v^{0,j}_sv^{1,j}_s$, with $j=0,1$. It is straighforward to see that the resulting graph is a hamilton cycle in $G$.

We now consider the case $\lambda \ge 1$ and $\mu> 1$. 
The list of paths involved in the construction is the following:

\begin{itemize}
\item in $K_0$ we take the path  $P$ in $H_{0,0}$ and in $H_{\mu,0}$;
we consider the subpths $v_s\,P\,v_h$ and $u_z\,P\,u_0$ in $H_{i,0}$ with $1\le i\le\mu-1$;

\item in $K_{\lambda}$ we consider  the path  $P$ in $H_{0,\lambda}$ and in $H_{\mu,\lambda}$;
we take the subpths $v_s\,P\, u_z$ and $v_k\,P\,u_0$ in $H_{i,\lambda}$ with $1\le i\le\mu-1$;
\end{itemize}

\noindent for $\lambda>1$,  we also consider the following subpaths:

\begin{itemize}
\item in $K_j$, with $j=1,2$, we take the path $P$ in $H_{0,j}$; the paths $v_s\,P\,v_h$ and $v_k\,P\,u_0$ in
$H_{1,j}$,  together with the isolated vertex $u^{1,j}_z$; we consider the paths  $v_s\,P\,v_h$ and $v_k\,P\,v_0$ in $H_{i,j}$, with $2\le i\le\mu-1$,
together with the isolated vertices $u^{i,j}_0$, $u^{i,j}_z$; we take the path $v_s\,P\,v_0$ in $H_{\mu,j}$ together with the isolated vertex $u^{\mu,j}_0$;

\item in $K_j$, with $3\le j\le\lambda-1$ and $\lambda>3$, we consider the path $v_s\,P\,u_1$ and the edge $u_0v_0$ in $H_{0,j}$; 
the paths  $v_s\,P\,v_h$ and $v_k\,P\,v_0$ in $H_{i,j}$, with $1 \le i\le\mu-1$,
together with the isolated vertices $u^{i,j}_0$, $u^{i,j}_z$; we take the path $v_s\,P\,v_0$ in $H_{\mu,j}$ together with the isolated vertex $u^{\mu,j}_0$.
\end{itemize}

In order to find a hamilton cycle in $G$, we connect the above subpaths by adding vertical and horizontal edges as we are going to describe.  We point out that we can join the subpaths by adding edges connecting the inner (outer) vertices in $P_{i, j}$ 
to the inner (outer) vertices in $P_{i+1,j}$ (respectively,   $P_{i,j+1}$) that are assigned the same label.
For instance,  we can join the subpath $v_s\,P_{i,j}v_h$ to the subpath $v_s\,P_{i+1,j}v_h$ by adding the edge with both end-vertices equal to $v_s$ (or $v_h$), 
because the vertex $v_s$ in $P_{i,j}$ corresponds to the vertex $v^{i,j}_s=v_{s+ib+ja}$, and the vertex $v_s$ in $P_{i+1,j}$ corresponds to the vertex $v^{i+1,j}_s = v_{s+(i+1)b+ja}$, 
and the two vertices are adjacent in $G$.  Analogously, for the outer vertices in $P_{i, j}$  and $P_{i, j+1}$.  We add the following edges:

\begin{itemize}

\item the vertical edges from $v_s$ in $H_{i,j}$ to $v_s$ in $H_{i+1, j}$ for every $0\le j\le\lambda$ and $0\le i\le\mu-1$, $i$ even;

\item the vertical edges from $v_h$ in $H_{i,j}$ to $v_h$ in $H_{i+1, j}$ for every $0\le j\le\lambda-1$ and $0\le i\le\mu-2$, $i$ odd;

\item  the vertical edges from $v_0$ in $H_{i,j}$ to $v_0$ in $H_{i+1, j}$ for $j=1,2$ and $2\le i\le\mu-1$, $i$ even, and  for $\lambda >3$  also
for every $3\le j\le\lambda-1$ and  $0\le i\le\mu-1$, $i$ even;

\item the vertical edges from $v_k$ in $H_{i,j}$ to $v_k$ in $H_{i+1, j}$ for every $1\le j\le\lambda$ and  $1\le i\le\mu-2$, $i$ odd;

\item the horizontal edges from $u_0$ in $H_{i,j}$ to $u_0$ in $H_{i, j+1}$ for every $2\le i\le\mu$ and $0\le j\le\lambda-1$,
and also for $i=0$ and $0\le j\le\lambda-1$, $j$ even, and for $i=1$ and $2\le j\le\lambda-1$,  together with the edge $u^{1,0}_0u^{1, 1}_0$;


\item  the horizontal edges from $u_z$ in $H_{i,j}$ to $u_z$ in $H_{i, j+1}$ for every $1\le i\le\mu-1$ and $0\le j\le\lambda-1$;

\item the horizontal edges from $u_1$ in $H_{0,j}$ to $u_1$ in $H_{0, j+1}$ for every  $3\le j\le\lambda-2$, $j$ odd  and $\lambda >3$.

\end{itemize}

\begin{figure}[htb]
\begin{center}
\includegraphics[width=15cm]{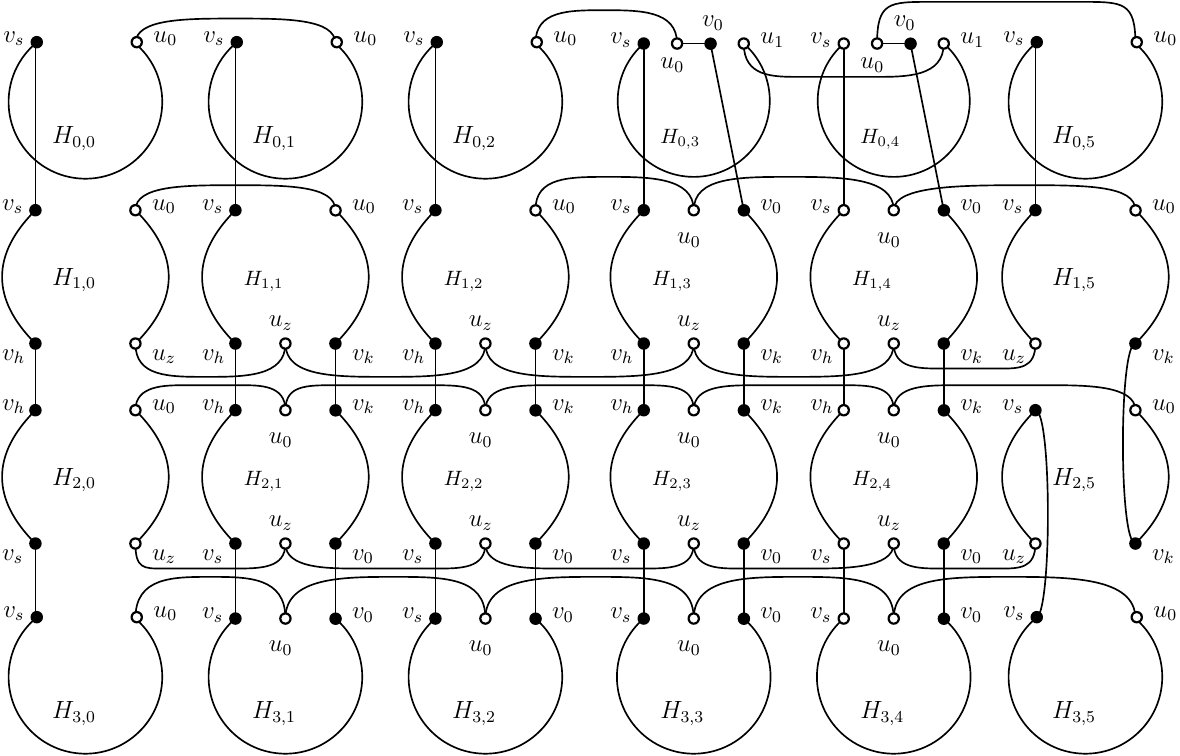}
\end{center}
\caption{The construction of a hamilton cycle in $G$ in Lemma  \ref{pro_method1_case}, when $\lambda$, $\mu$ are both odd; 
in the figure  we have $\lambda=5$ and $\mu=3$. }\label{fig_ham_cycle_lambda_mu_odd}
\end{figure}

We show the construction for $\lambda=5$, $\mu=3$ in Figure \ref{fig_ham_cycle_lambda_mu_odd}.
We prove that it yields a hamilton cycle in $G$ when $\lambda>1$; the proof for the case $\lambda=1$ is 
similar and will be omitted.  We denote with $U$ the subset of $V(G)$ consisting of the copies of the vertices 
$u_0$, $u_z$ in each $H_{i,j}$,  i.e., $U=\{u^{i,j}_0, u^{i,j}_z : 0\le i\le\mu, 0\le j\le\lambda\}$; 
let $U'$ be the subset of $U$ defined by $U'=\{u^{0,1}_0, u^{1,1}_0,  u^{1,0}_0,  u^{0,j}_0: 3\le j\le\lambda\}\cup\{u^{1,j}_z :    0\le j\le\lambda\}$.  

We consider the paths in the vertical components $K_0$, $K_2$, $K_{\lambda}$  of $G$ that are different from the paths $u_0\,P_{1,0}\,u_z$, $P_{0,\lambda}$, $v_s\,P_{1,\lambda}u_z$, together with the vertices in $U\setminus U'$, 
and the above horizontal and vertical edges joining them. They form a path from $u^{0,0}_0$ to $u^{0,2}_0$ covering all vertices in  the union of the vertical components $K_0$, $K_2$, $K_{\lambda}$  together with the vertices in  
$U\setminus U'$,
with the exception for the vertices belonging to the paths $u_0\,P_{1,0}\,u_z$, $P_{0,\lambda}$, $v_s\,P_{1,\lambda}u_z$; we denote such a path with $u^{0,0}_0\,P'\,u^{0,2}_0$.

Now we consider the paths  $u_0\,P_{1,0}\,u_z$, $P_{0,\lambda}$, $v_s\,P_{1,\lambda}u_z$, together with the vertices in $U'$, the paths we defined for the vertical components $K_1$ and $K_j$, with $3\le j\le\lambda-1$ with $\lambda>3$,   
and the above horizontal and vertical edges joining them. They form a path from $u^{0,1}_0$ to $u^{0,3}_0$ covering all vertices belonging to their union; we denote such a path with $u^{0,1}_0\,P'\,u^{0,3}_0$.
We join the paths $u^{0,0}_0\,P'\,u^{0,2}_0$,  $u^{0,1}_0\,P'\,u^{0,3}_0$ trough the  edges $u^{0,0}_0u^{0,1}_0$, $u^{0,2}_0u^{0,3}_0$, 
and form a hamilton cycle in $G$.

\end{proof}

We now introduce the second construction of this section.

\begin{lemma}\label{pro_method1_case2}
Let $G=B(m; a, S, b)$ be a connected bicirculant graph with  
$|S|\geq 4$ and 
$a,b\ne m/2$  
that has a $(\lambda, \mu)$-uniform representation such that $\lambda=\lambda(b)=\gcd(m, S, b)-1 >0$ and
$\mu=\mu(b)=\gcd(m, S)/\gcd(m, S, b)-1\ge 0$.
Assume that $\lambda$ is even and $\mu$ is odd, or $\lambda$ is odd and 
$\mu$ is even. Then $G$ is hamiltonian.
\end{lemma}

\begin{proof}
If the connected components of $H(m;S)$ are hamiltonian, then the assertion follows from Lemma \ref{lem_BoPiZi2025_2ndpaper}. Now assume that this is not the case. The uniform representation of $G$ allows us to describe it as detailed at the beginning of this section and, in what follows, we will use this representation to construct a hamilton cycle in $G$. By Lemma \ref{lemma_path_1} we can assume that $H_{0,0}$ contains a hamilton path $P$  described earlier in this section.

We consider the cases where $\lambda$ is even and $\mu$ is odd or $\lambda$ is odd and $\mu=0$.
We omit the case where $\lambda$ is odd and $\mu$  is even, $\mu>0$, because it is equivalent to the case where $\lambda$ is even, $\lambda>0$, and  $\mu$ is odd. In fact, it suffices to replace the outer vertices with the inner vertices, and vice versa, in the construction that we are going to define; see the explanation given in the proof of Lemma \ref{pro_method1_case}.
The case where $\lambda$ is odd and $\mu=0$ cannot be obtained from the case where $\lambda$ is even and $\mu$ is odd by exploiting the symmetry between outer and inner vertices, since we are considering $\lambda>0$. This case -- $\lambda$ odd and $\mu=0$ -- is studied in Case 3.  
For $\lambda$ even we treat separately the cases $\mu=1$, and $\mu>1$.

For each of the cases, we will  list the paths that are used in the construction and specify the edges joining pairs of paths.
We omit showing that they form a hamilton cycle
in $G$, since the proof is similar to the proof  of Lemma \ref{pro_method1_case} for the case  $\lambda>1$.
\medskip

\noindent
{\bf Case 1:  $\lambda$ is even and $\mu=1$.}
We recall that the vertices $u_z$ and $v_s$ in the definition of the path $P$ are adjacent.  We will use the edge $u_zv_s$ to define the hamilton path $v_h\,P'\,u_0=v_h\,P\,v_s\cup v_su_z\cup u_z\,P\,u_0$
from $v_h$ to $u_0$, and the subpath $v_s\,P''\,u_0=v_su_z\cup u_z\,P\,u_0$ that spans $G$ together with the subpath $u_t\,P\,v_h$.  We show the construction for $\lambda=4$, $\mu=1$ in Figure \ref{fig_ham_cycle_lambda_even_mu=1}, which can be used by the reader to check that the paths involved in the construction of a hamilton cycle in $G$ when $\mu=1$ are the following:

\begin{itemize}
\item in $K_0$ we take the path $v_h\,P'\,u_0$ 
in $H_{0,0}$ and $H_{1,0}$;

\item in $K_j$, with $1\le j\le\lambda$,  $j\neq \lambda-1$, we consider the subpaths $v_s\,P\,v_h$, $u_z\,P\,u_0$ in $H_{0,j}$,
and  $v_s\,P''\,u_0$, 
$u_t\,P\,v_h$ in $H_{1,j}$;

\item in $K_{\lambda-1}$ we consider $v_s\,P\,v_h$, $u_z\,P\,v_0$ in $H_{0,\lambda-1}$,  the paths
$u_t\,P\,v_h$ and $v_s\,P''\,v_0=v_su_z\cup u_z\,P\,v_0$  in $H_{1,\lambda-1}$,  together with the isolated vertex $u_0$ in $H_{0,\lambda-1}$ and $H_{1,\lambda-1}$.

\end{itemize}

\noindent We add the following edges to form a hamilton cycle in $G$:
$v^{0,j}_sv^{1,j}_s$ for $1\le j\le\lambda$, and $v^{0,j}_hv^{1,j}_h$ for $0\le j\le\lambda$;
 the edeges $v^{0,\lambda-1}_0v^{1,\lambda-1}_0$;
 the edges $u^{i,j}_0u^{i,j+1}_0$ for $i=0,1$, $0\le j\le\lambda-2$, $j$ even, and $u^{i,\lambda-1}_0u^{i,\lambda}_0$;
 the edges $u^{0,j}_zu^{0,j+1}_z$ and $u^{1,j}_tu^{1,j+1}_t$ for $1\le j\le\lambda-1$, $j$ odd.
\smallskip

\begin{figure}[t]
\begin{center}
\includegraphics[width=12cm]{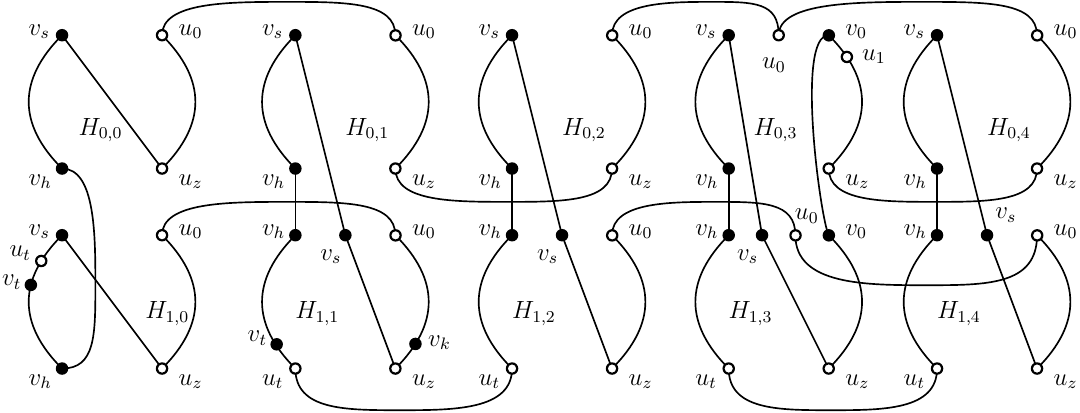}
\end{center}
\caption{The construction of a hamilton cycle in $G$ described in Lemma  \ref{pro_method1_case2} when $\lambda$ is even and $\mu=1$; in the figure we have $\lambda=4$.}\label{fig_ham_cycle_lambda_even_mu=1}
\end{figure}

\noindent
{\bf Case 2:  $\lambda$ is even and $\mu$ is odd, $\mu>1$.}
In Figure \ref{fig_ham_cycle_lambda_even_mu_odd} we show the construction for $\lambda=4$ and $\mu=3$.  
The list of paths involved in the construction is the following:

\begin{itemize}
\item in $K_0$ we take the path  $P$ in $H_{0,0}$ and in $H_{\mu,0}$;
in $H_{i,0}$, with $1\le i\le\mu-1$, we consider the subpths $v_s\,P\,v_h$ and $u_z\,P\,u_0$;

\item in each $K_j$, with $1\le j\le\lambda-2$,   in $H_{0,j}$ we take the edge $u_0v_0$ and the subpath $v_s\,P\,u_1$;
in  $H_{i,j}$, with $1\le i\le\mu-3$, we consider  the subpaths $v_s\,P\,v_h$ and $u_z\,P\,u_1$  together with the isolated vertices $u^{i,j}_0$ and $v^{i,j}_0$;
in $H_{\mu-2,j}$ we take the subpaths $v_s\,P\,v_h$, $v_k\,P\,u_1$ and the isolated vertices $u^{\mu-2,j}_0$, $v^{\mu-2,j}_0$, $u^{\mu-2,j}_z$;
in $H_{\mu-1,j}$, we remove the edges incident to the copies of $u_0$ and $u_z$ and consider the remaining subpaths $v_s\,P\,v_h$,  $v_k\,P\,v_0$ together with the
isolated vertices $u^{\mu-1,j}_0$, $u^{\mu-1,j}_z$; we take the path $P$ in $H_{\mu,j}$;

\item in $K_{\lambda-1}$ we consider the following subpaths: $v_s\,P\,v_0$ and the isolated vertex $u^{i,\lambda-1}_0$ in $H_{i,\lambda-1}$ with $i=0,\mu$; the paths 
$v_s\,P\,v_h$, $v_k\,P\,v_0$ and the isolated vertices $u^{i,\lambda-1}_0$, $u^{i,\lambda-1}_z$ in $H_{i,\lambda-1}$ with $1\le i\le\mu-3$; the paths
$v_s\,P\,v_h$, $u_z\,P\,v_0$ and the isolated vertex $u^{i,\lambda-1}_0$ in  $H_{i,\lambda-1}$ with $i=\mu-2, \mu-1$;

\item in $K_{\lambda}$ we take the path $P$ in $H_{i,\lambda}$ with $i=0,\mu-2, \mu-1,\mu$ and  the subpaths $v_s\,P\,u_z$, $v_k\,P\,u_0$ in $H_{i,\lambda}$ with $1\le i\le\mu-3$.
\end{itemize}

In order to find a hamilton cycle in $G$, we connect the above subpaths by adding vertical and horizontal edges as follows:

\begin{itemize}

\item the vertical edges from $v_s$ in $H_{i,j}$ to $v_s$ in $H_{i+1, j}$ for every $0\le j\le\lambda$ and $0\le i\le\mu-1$, $i$ even;

\item the vertical edges from $v_h$ in $H_{i,j}$ to $v_h$ in $H_{i+1, j}$ for every $0\le j\le\lambda-1$ and $1\le i\le\mu-2$, $i$ odd;

\item  the vertical edges from $v_0$ in $H_{i,j}$ to $v_0$ in $H_{i+1, j}$ for every $1\le j\le\lambda-2$ and $0\le i\le\mu-2$,
and for $j=\lambda-1$ and $0\le i\le\mu-1$, $i$ even;

\item the vertical edges from $v_k$ in  in $H_{i,j}$ to $v_k$ in $H_{i+1, j}$ for every $1\le j\le\lambda-2$ and $i=\mu-2$, and 
for $\lambda-1\le j\le\lambda$ and  $1\le i\le\mu-4$, $i$ odd;

\item the horizontal edges from $u_0$ in $H_{i,j}$ to $u_0$ in $H_{i, j+1}$ for every $1\le i\le\mu-1$ and $0\le j\le\lambda-1$,
and also for $i=0, \mu$ and $0\le j\le\lambda-2$, $j$ even, together with $j=\lambda-1$;

\item  the horizontal edges from $u_z$ in $H_{i,j}$ to $u_z$ in $H_{i, j+1}$ for every $1\le i\le\mu-3$ and $0\le j\le\lambda-2$, $j$ even, together with $j=\lambda-1$, 
and also for $\mu-2\le i\le\mu-1$ and $0\le j\le\lambda-2$;

\item the horizontal edges from $u_1$ in $H_{i,j}$ to $u_1$ in $H_{i, j+1}$ for every $0\le i\le\mu-2$ and $1\le j\le\lambda-3$, $j$ odd.

\end{itemize}
\smallskip

\begin{figure}[htb]
\begin{center}
\includegraphics[width=12cm]{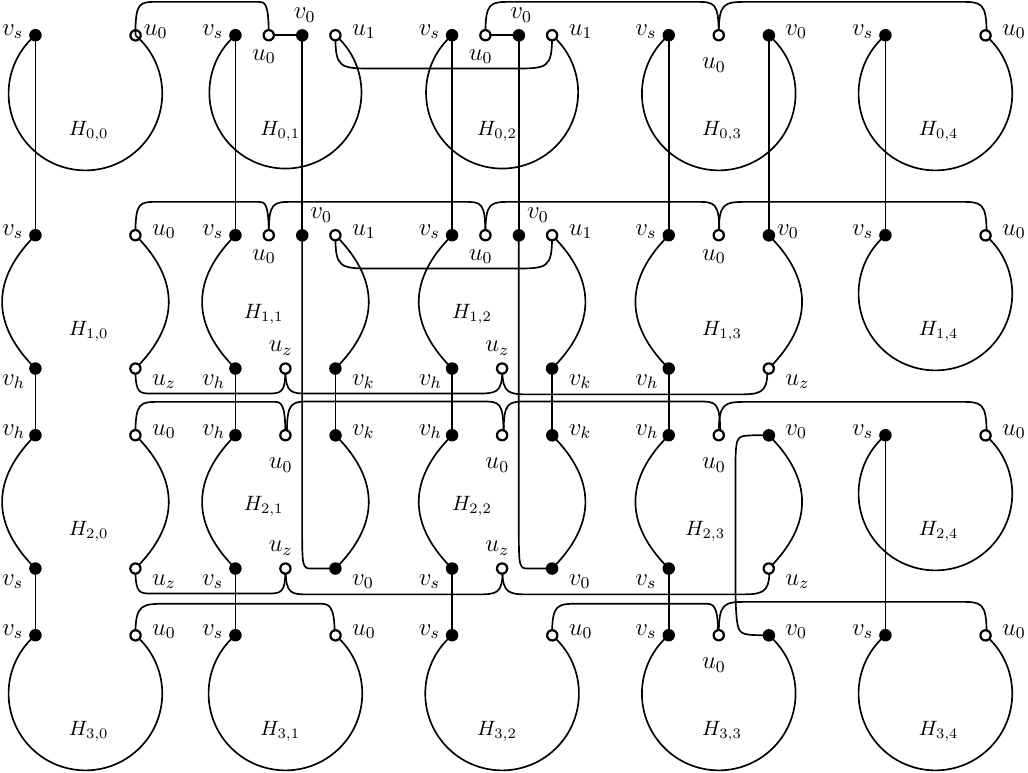}
\end{center}
\caption{The construction of a hamilton cycle in $G$ described in Lemma  \ref{pro_method1_case2} when $\lambda$ is even and $\mu$ is odd, $\mu>1$; in the figure  we have $\lambda=4$ and $\mu=3$.}\label{fig_ham_cycle_lambda_even_mu_odd}
\end{figure}

\noindent
{\bf Case 3: $\lambda$ is odd and $\mu=0$.}
Since $\mu=0$ we have that $\gcd(m, S)=\gcd(m, S, b)>1$. Each component $K_j$ of $K$ contains only the component $H_{0,j}$ of $H$, for $0\leq j\leq\lambda$. This means that the path $P$  in $K_0 $ (or $H_{0,0}$) has chords of type $b$.  We can repeat the construction given in the proof of Lemma \ref{lem_BoPiZi2025_2ndpaper} (see \cite[Lemma 3.2]{BoPiZi2025_2ndpaper}).
For completeness, we report it below.

We will denote with $u_r$, $u_{\ell}$ and $u_p$, $u_q$ the vertices adjacent to $v_{s-b}$ and $v_{s+b}$, and assume that they are arranged on 
$P$ in the order $ v_s$,  $u_t$,  $u_r$, $v_{s-b}$,  $u_{\ell}$,  $u_p$, $v_{s+b}$, $u_q$, $u_0$.

In this setting we can find a hamilton path $u_0\,P'\,u_p$ in $H_{0,0}$ from $u_0$ to $u_p$ that is defined as follows:  $u_0\,P'\,u_p=u_0\,P\,v_{s+b}\cup\, v_{s+b}v_s\,\cup v_s\,P\,u_p$. We will also consider the path 
$u_0\,P''\,u_t$ from $u_0$ to $u_t$ given by $u_0\,P''\,u_t=u_0\,P\,v_{s+b}\cup (v_{s+b}, v_s,  v_{s-b})\cup v_{s-b}\,P\,u_t$; 
the union of the paths  $u_p\,P\,u_{\ell}$ and   $u_0\,P''\,u_t$  spans the vertices of $K_0$ (or $H_{0,0}$).

We use the above paths to construct a hamilton cycle in $G$.  More specifically, in $H_{0,0}$ and $H_{0,\lambda}$ we consider the copy of
$u_0\,P'\,u_p$;
for $1\le j\le\lambda-1$,  we consider the paths $u_p\,P_{0,j}\,u_{\ell}$ and $u_0\,P''_{0,j}\,u_t$ in $H_{0,j}$. We find a hamilton cycle in $G$ by joining the paths through the following edges: the edges from $u_0$, $u_p$ in $H_{0,j}$ to $u_0$, $u_p$ in $H_{0,j+1}$, respectively, for $0\le j\le\lambda-1$, $j$ even, and the edges from $u_t$, $u_{\ell}$ in $H_{0,j}$ to $u_t$, $u_{\ell}$ in $H_{0,j+1}$, respectively, for $1\le j\le\lambda-2$, $j$ odd.
In particular, for $\lambda=1$, it suffices to join the hamilton path $u_0\,P'\,u_p$ in $H_{0,0}$ to the hamilton path $u_0\,P'\,u_p$ in $H_{0,1}$ through the edges from $u_x$ in $H_{0,0}$ to $u_x$ in $H_{0,1}$, with $x \in \{0, p\}$. 

\end{proof}

Now we introduce the third construction of this section. This construction will be needed only for odd $m$. When $m$ is odd, $\lambda$ and $\mu$ are both even. The construction, however, works only when $\mu>0$. This is why we require the assumption $\gcd(m,S,b)<\gcd(m,S)$ in Theorem \ref{th_main} for odd $m$.

\begin{lemma}\label{pro_method1_case3}
Let $G=B(m; a, S, b)$ be a connected bicirculant graph with 
$|S|\ge 4$ and  
$a, b\ne m/2$ 
that has a $(\lambda, \mu)$-uniform representation such that $\lambda=\lambda(b)=\gcd(m, S, b)-1 >0$ and
$\mu=\mu(b)=\gcd(m, S)/\gcd(m, S, b)-1\ge 0$.
Assume that $\lambda$ and $\mu$ are even and 
$\mu > 0$.
Then $G$ is hamiltonian.
\end{lemma}

\begin{proof} 
If the connected components of $H(m;S)$ are hamiltonian, then the assertion follows from Lemma \ref{lem_BoPiZi2025_2ndpaper}. Now assume that this is not the case. The uniform representation of $G$ allows us to describe it as detailed at the beginning of this section and, in what follows, we will use this representation to construct a hamilton cycle in $G$. By Lemma \ref{lemma_path_1} we can assume that $H_{0,0}$ contains a hamilton path $P$  described earlier in this section. 

We give the list of paths that are involved in the construction in the case where $\lambda >2$ and $\mu\ge 2$; the same list can be used even in the case where $\lambda=\mu=2$. The case with $\lambda=2$ and $\mu>2$ can be obtained from the case with $\lambda>2$ and $\mu=2$ by exchanging the inner vertices with outer vertices; see the explanation given in the proof of Lemma \ref{pro_method1_case}.
The list of paths is the following:

\begin{itemize}
\item in $K_0$:  \begin{itemize} 
\item the hamilton path $v_h\,P'\,u_0=v_h\,P\,v_s \cup\, v_su_z\,\cup u_z\,P\,u_0$ in $H_{0,0}$;
\item the subpaths $v_s\,P\,v_h$ and $u_z\,P\,u_0$ in $H_{i,0}$ with $1\le i\le\mu-1$;
\item  $P$ in $H_{\mu,0}$;
\end{itemize}

\item in $K_j$, with $1\le j\le\lambda-2$,  $\lambda>2$: 
\begin{itemize} 
\item  $P$ in $H_{0,j}$;
\item the subpaths $u_t\,P\,u_z$ and $v_k\,P\,v_0$ in $H_{\mu-1,j}$ togheter with the isolated vertices $u_0$, $v_s$; 
\item the subpaths $v_s\,P\,v_k$ and $u_k\,P\,v_0$ together with the isolated vertex $u_0$ in $H_{\mu,j}$;

\item if $\mu>2$, then for $j=1$ we also consider the supbaths  $v_s\,P\,v_h$ and $u_z\,P\,v_0$ together with the isolated vertex $u_0$ in each $H_{i,1}$ with $1\le i\le\mu-2$;
for $2\le j\le\lambda -2$, we also consider the supbath  $v_s\,P\,v_0$ together with the isolated vertex $u_0$ in each $H_{i,j}$ with $1\le i\le\mu-2$;
\end{itemize}

\item in $K_{\lambda-1}$:
\begin{itemize} 
\item the subpath $v_t\,P\,v_0$, the edge $u_tv_s$ and the isolated vertex $u_0$ in $H_{i,\lambda-1}$ for every $0\le i\le\mu-2$, $i$ even;

\item the edges $u_xv_x$ for $1\leq x\le k$ and $x=t$, the subpath  $u_{t-1}\,P\,u_z$, and the isolated vertices $v_s$, $u_0$, $v_0$ in $H_{\mu-1,\lambda-1}$;

\item the edges $v_xu_{x+1}$ for $0\leq x\le k-1$,  the subpath  $v_s\,P\,v_k$, and the isolated vertex $u_0$ in $H_{\mu,\lambda-1}$;

\item the subpaths $u_{t-1}\,P\,v_0$,  $(v_t, u_t, v_s)$ and the isolated vertex $u_0$ in $H_{i,\lambda-1}$ for every $1\le i\le\mu-3$, $i$ odd and $\mu>2$;

\end{itemize}

\item in $K_{\lambda}$:
\begin{itemize} 
\item  the paths $(u_t, v_s, u_z)\cup u_z\,P\,u_0$ and $v_t\,P\,v_h$ in $H_{i,\lambda}$ for every $0\le i\le\mu-2$, $i$ even;
\item the edges $u_xv_x$ for $0\leq x\le k-1$ and $x=t$,  the path $(u_k, v_k, u_z, v_s)$,  and the subpath  $u_{t-1}\,P\,v_h$ in $H_{\mu-1,\lambda}$;
\item the edges $u_xv_x$ for $0\leq x\le k-1$,  and the subpath  $v_s\,P\,u_k$ in $H_{\mu,\lambda}$.
\item the paths $(v_t,  u_t,  v_s,  u_z)\cup u_z\,P\,u_0$ and $u_{t-1}\,P\,v_h$ in $H_{i,\lambda}$ for every $1\le i\le\mu-3$, $i$ odd and $\mu>2$;
\end{itemize}
\end{itemize}

We add the following edges in order to form a hamilton cycle in $G$:

\begin{itemize}

\item  the vertical edges from $v_s$ in $H_{i,j}$ to $v_s$ in $H_{i+1, j}$ 
for $j=0$ and  $1\le i\le \mu-1$, $i$ odd, for $1\le j\le\lambda-1$ and $0\le i\le\mu-2$, $i$ even, 
for $1\le j\le\lambda$ and $i=\mu-1$;


\item the vertical edges from $v_h$ in $H_{i,j}$ to $v_h$ in $H_{i+1, j}$ for $j\in\{0, \lambda\}$ and $0\le i\le\mu-2$, $i$ even, and also for
$j=1$ and  $1\le i\le\mu-3$, $i$ odd and $\mu>2$;

\item  the vertical edges from $v_0$ in $H_{i,j}$ to $v_0$ in $H_{i+1, j}$ for every $1\le j\le\lambda-2$ with $1\le i\le\mu-1$, $i$ odd, and for $j=\lambda-1$ with $0\le i\le\mu-2$, $i$ even;

\item the vertical edges from $v_k$ in $H_{\mu-1,j}$ to $v_k$ in $H_{\mu, j}$ for every $1\le j\le\lambda-1$;

\item the vertical edges from $v_t$ in $H_{i,j}$ to $v_t$ in $H_{i+1, j}$ for $j=\lambda-1, \lambda$ and $0\le i\le\mu-2$, $i$ even;

\item the vertical edges from $v_x$ in $H_{\mu-1,j}$ to $v_x$ in $H_{\mu, j}$ for $j=\lambda-1, \lambda$ and $0\le x\le k-1$;

\item the horizontal edges from $u_0$ in $H_{i,j}$ to $u_0$ in $H_{i, j+1}$ for every $1\le i\le\mu$ and $0\le j\le\lambda-1$,
for $i=0$ and $0\le j\le\lambda-2$, $j$ even, together with $i=0$ and $j=\lambda-1$;


\item  the horizontal edges from $u_z$ in $H_{i,j}$ to $u_z$ in $H_{i, j+1}$ for $i=\mu-1$ and $0\le j\le\lambda-2$, $j$ even, for every $1\le i\le\mu-2$  and $j=0$ if $\mu>2$;

\item the horizontal edges from $u_t$ in $H_{i,j}$ to $u_t$ in $H_{i, j+1}$ for every $0\le i\le\mu-2$, $i$ even and $j=\lambda-1$,  and also for 
$i=\mu-1$ and $1\le j\le\lambda-1$, $j$ odd;

\item the horizontal edges from $u_{t-1}$ in $H_{i,\lambda-1}$ to $u_{t-1}$ in $H_{i, \lambda}$ for every $1\le i\le\mu-1$, $i$ odd;

\item the horizontal edges from $u_{x}$ in $H_{i,\lambda-1}$ to $u_x$ in $H_{i, \lambda}$ for $i=\mu-1, \mu$ and $1\le x\le k$;

\item the horizontal edges from $u_{k}$ in $H_{\mu,j}$ to $u_k$ in $H_{\mu, j+1}$ for $1\le j\le\lambda-3$, $j$ odd.

\end{itemize}

Figure \ref{fig_ham_cycle_even} shows the paths for the case $\lambda=\mu=4$; the path $u_z\,P\,u_0$ contains $7$ vertices ($k=2$), 
the path $v_s\,P\,u_z$ has at least $6$ vertices. 

\begin{figure}[h]
\begin{center}
\includegraphics[width=13cm]{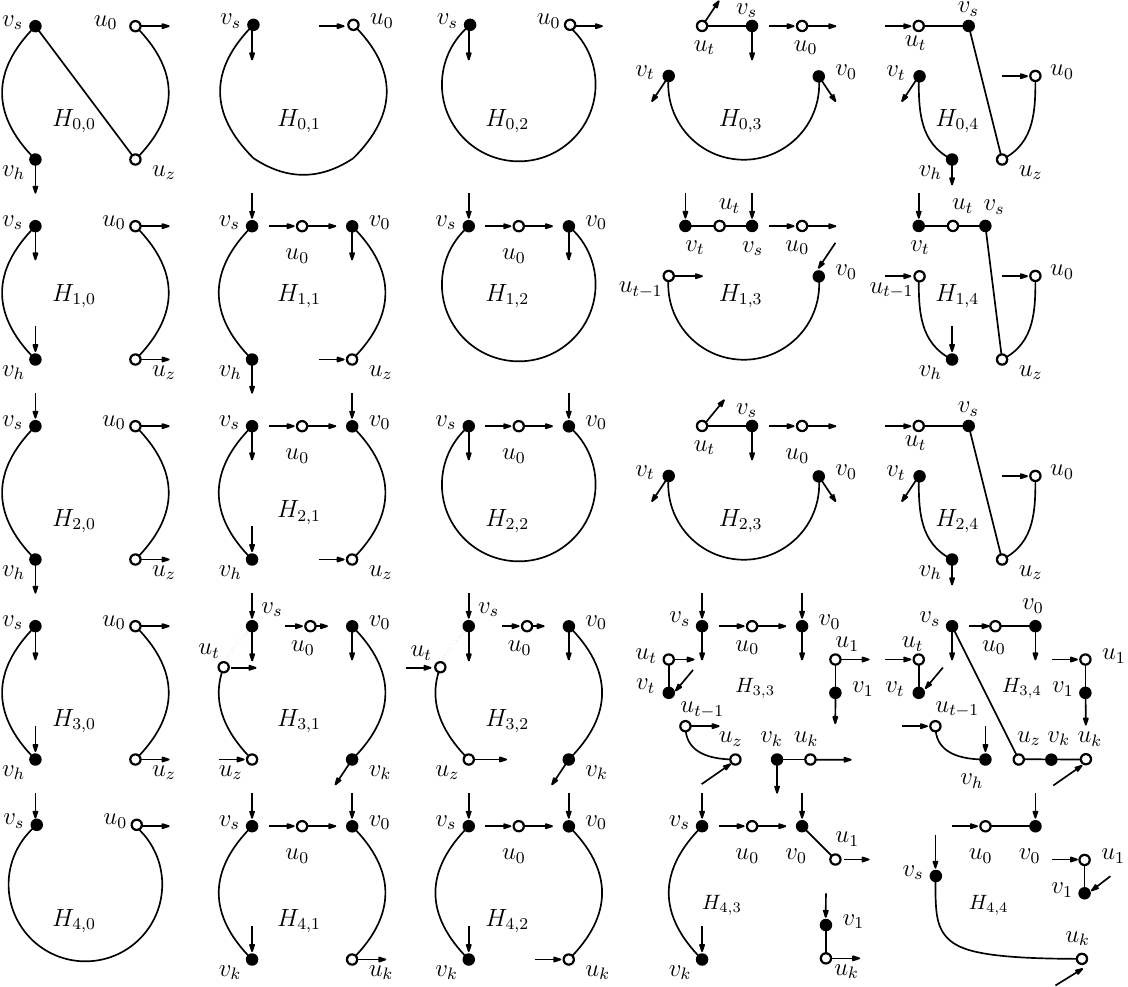}
\end{center}
\caption{The construction of a hamilton cycle in  Lemma  \ref{pro_method1_case3},  when $\lambda$ and $\mu$ are both even, $\lambda>2$ and $\mu\ge 2$. In the figure we have $\lambda=\mu=4$; the path $u_z\,P\,u_0$ of the figure contains $7$ vertices,  that is $k=2$, 
the path $v_s\,P\,u_z$ has at least $6$ vertices. The arrows represent the edges connecting the vertices that belong to different copies of $H_{0,0}$:
an arrow leaving (entering) $u_x$ in $H_{i ,j}$ means that the copy of $u_x$ in $H_{i,j}$ is connected to the copy of $u_x$ in $H_{i,j+1}$ (respectively, in $H_{i,j-1}$); an arrow leaving (entering) $v_y$ in $H_{i ,j}$ means that the copy of $v_y$ in $H_{i,j}$ is  connected to the copy of $v_x$ in $H_{i+1,j}$ (respectively, in $H_{i-1,j}$). 
}
\label{fig_ham_cycle_even}
\end{figure}

The rest of the proof is similar to that of Lemma \ref{pro_method1_case} and can be adapted to the case $\lambda>2$
under consideration, as we are going to explain. More specifically, the vertical edges we add as described above produce the following paths. In $K_0$, we obtain a path $u^{0,0}_0\,P'\,u^{\mu, 0}_0$ from $u^{0,0}_0$ to $u^{\mu, 0}_0$ that covers all vertices of $K_0$ with the exception of the copy of the path $u_z\,P\,u_0$ in each $H_{i,0}$ with $1\le i\le\mu-1$. 
In $K_1$, we obtain a path $u^{0,1}_0\,P'\,u^{\mu,1}_k$ 
from $u^{0,1}_0$ to $u^{\mu,1}_k$ that covers all vertices of $K_1$ with the exception of the vertices $u^{i,1}_0$ with $1\le i\le\mu$, of the copy of the path $u_t\,P\,u_z$ in $H_{\mu-1,1}$, and of the copy of the path $u_z\,P\,v_0$ in each $H_{i, 1}$ with $1\le i\le\mu-2$ if $\mu>2$.  In $K_j$, with $2\le j\le\lambda-2$ ($\lambda>2$), we obtain the path $u^{0,j}_0\,P'\,u^{\mu,j}_k$ from the vertex $u^{0,j}_0$ to the vertex $u^{\mu,j}_k$ that covers all vertices in $K_j$ with the exception of the copy of $u_0$ in $H_{i,j}$, with $1\le i\le\mu$, and of the copy of the path $u_t\,P\,u_z$ in $H_{\mu-1},j$. 
In $K_{\lambda-1}\cup K_{\lambda}$, the above vertical edges together with the horizontal edges from $u_t, u_{t-1},u_1,\dots, u_k$ in $H_{i,\lambda-1}$ to 
$u_t, u_{t-1},u_1,\dots, u_k$ in $H_{i,\lambda}$ allow us to find a path $u^{\mu-2,\lambda}_0\,P'\,u^{\mu-1,\lambda-1}_z$ from $u^{\mu-2,\lambda}_0$ to $u^{\mu-1,\lambda-1}_z$, and for $\mu\ge 4$ also a path $u^{i,\lambda}_0\,P'\,u^{i+1,\lambda}_0$ 
from the vertex $u^{i,\lambda}_0$ to the vertex $u^{i+1,\lambda}_0$ with $0\le i\le\mu-4$, $i$ even. The union of these paths covers all vertices in $K_{\lambda-1}\cup K_{\lambda}$ with the exception of the copy of $u_0$ in $H_{i,\lambda-1}$, $0\le i\le\mu$, and of the copy of the edge $u_0v_0$ in $H_{i,\lambda}$, $i=\mu-1, \mu$.

The addition of the edges from $u_0, u_z, u_t$ in $H_{i,j}$ to $u_0, u_z, u_t$ in $H_{i,j+1}$, as described above,
connects the paths  $u^{0,0}_0\,P'\,u^{\mu, 0}_0$ and  $u^{0,1}_0\,P'\,u^{0,1}_k$
to the paths in $K_{\lambda-1}\cup K_{\lambda}$ and forms a path from $u^{\mu, 1}_k$ to  $u^{\mu,\lambda-2}_k$, 
denoted as $u^{\mu, 1}_k\,\overline P\,u^{\mu,\lambda-2}_k$. It covers all vertices in $K_0\cup K_1\cup K_{\lambda-2}\cup K_{\lambda-1}\cup K_{\lambda}$, and for $\lambda> 4$ also covers the copy of the vertex $u_0$ in each $H_{i,j}$ with $1\le i\le\mu$ and $2\le j\le\lambda-3$. 
For $\lambda=4$, the addition of the edge $u^{\mu, 1}_k\,u^{\mu,\lambda-2}_k$ turns the path $u^{\mu, 1}_k\,\overline P\,u^{\mu,\lambda-2}_k$ into a hamilton cycle of $G$, and the assertion follows. For $\lambda>4$, the addition of the edges from $u_0, u_z, u_t$ in $H_{i,j}$ to $u_0, u_z, u_t$ in $H_{i,j+1}$ also connects the paths $u^{0,j}_0\,P'\,u^{\mu,j}_k$ with $2\le j\le\lambda-3$ and form a path $u^{\mu,j}_k\,\overline P\,u^{\mu,j+1}_k$ from the vertex $u^{\mu,j}_k$ to $u^{\mu,j+1}_k$ with $2\le j\le\lambda-4$, $j$ even. The addition of the horizontal edges from $u_k$ in $H_{\mu, j}$ to $u_k$ in $H_{\mu,j+1}$, with $1\le j\le\lambda-2$, $j$ odd, allows us to connect the paths $u^{\mu, 1}_k\,\overline P\,u^{\mu,\lambda-2}_k$ and $u^{\mu,j}_k\,\overline P\,u^{\mu,j+1}_k$,  with $2\le j\le\lambda-4$, $j$ even, and find a hamilton cycle in $G$. The assertion follows.
\end{proof}

\section{Assembling the main results}
\label{sec_further_constructions}

In this section, we prove our main results. We first prove Theorem \ref{th_main}. To this end, we combine the constructions of Lemmas 
\ref{lem_BoPiZi2025_2ndpaper}--\ref{pro_method1_case3} with the construction given in Lemma 3.3 from \cite{BoPiZi2025_2ndpaper}. This yields the construction described in the proof of Theorem \ref{th_main}. 
We stress that Theorem \ref{th_main} is a stronger version of \cite[Lemma 3.4]{BoPiZi2025_2ndpaper}. In particular, in \cite[Lemma 3.4]{BoPiZi2025_2ndpaper} we require  the connected components of $H(m; S)$ to be hamiltonian, whereas in Theorem \ref{th_main} we drop this assumption.
For completeness, we restate Lemma 3.3 from \cite{BoPiZi2025_2ndpaper}.

\begin{lemma}[\protect{\cite[Lemma 3.3] {BoPiZi2025_2ndpaper}}]\label{pro_construction_Bd1bis}
Let $G = B(m; R, S, T)$ be a connected bicirculant graph  of degree $d \geq 4$, with $|R|,|T|\geq 3$.
Let $H$ be the subgraph of $G$ obtained by removing the edges of types $a\in R$ and $b\in T$, with $a, b\neq m/2$.

If the connected components of $H$ contain a hamilton cycle with at least two outer edges, then $G$ is hamiltonian.
\end{lemma}

The following remark concerning the constructions described in \cite{BoPiZi2025_2ndpaper}
is in order. 
 
\begin{remark}\label{rem_construction}
The construction from Lemma \ref{lem_BoPiZi2025_2ndpaper} yields a hamilton cycle  with at least two outer edges, see \cite[Lemma 3.2]{BoPiZi2025_2ndpaper}. 
From the proof of  \cite[Lemma 3.3] {BoPiZi2025_2ndpaper}, it can be seen that the construction given in \cite[Lemma 3.4] {BoPiZi2025_2ndpaper} also yields a hamilton cycle  with at least two outer edges.
\end{remark}

We are now ready to prove Theorem \ref{th_main}.

\begin{proof}[Proof of Theorem \ref{th_main}]
Recall that $G=B(m; R, S, T)$ is a connected bicirculant graph with $|S|\ge 4$ and $\gcd(m, S)>1$. Furthermore,  $R\cup T$ contains an element that is not coprime to $\gcd(m, S)$  and, if $m$ is odd,  is not divisible by $\gcd(m, S)$. We need to prove that $G$ is hamiltonian.

If $m/2\in R, T$, then the assertion follows from Theorem \ref{th_m/2_2ndpaper}, while it follows from Theorem \ref{th_bicirculant_m/gcd_3rdpaper} if $m/\gcd(m, S)$ is even. In the following, we consider the case where $m/\gcd(m,S)$ is odd and $m/2\not\in R\cup T$. 
We proceed  similarly as in the proof of Lemma 3.4 in \cite{BoPiZi2025_2ndpaper}.

We construct a regular spanning subgraph $L$ of $G$ that contains $H$ as follows. Let $b_1$ be an element of $R \cup T$ that is not coprime to $\gcd(m,S)$ and, if $m$ is odd, not divisible by $\gcd(m,S)$. Without loss of generality we may assume that $b_1 \in T$.  Denote by $a_1$ an arbitraty element of $R$ and define $R_1=\{a_1,-a_1\}$, $T_1=\{b_1,-b_1\}$ and $L_1=B(m;R_1,S,T_1)$. If the graph $L_1$ is connected, that is, $\gcd(m,R_1,S,T_1)=1$,  we set $L=L_1$. Otherwise pick arbitrary elements from $R \setminus \{a_1,-a_1\}$ and  $T \setminus \{b_1,-b_1\}$ and denote them by $a_2$ and $b_2$, respectively. Define $R_2=R_1 \cup\{a_2,-a_2\}$, $T_2=T_1 \cup\{b_2,-b_2\}$ and $L_2=B(m;R_2,S,T_2)$. If the graph $L_2$ is connected, that is, $\gcd(m,R_2,S,T_2)=1$, we set $L=L_2$. Otherwise we continue the process: at each step we pick arbitrary elements $a_{i+1} \in R\setminus R_i$ and  $b_{i+1} \in T\setminus T_i$ and define $R_{i+1}=R_i \cup\{a_{i+1},-a_{i+1}\}$, $T_{i+1}=T_i \cup\{b_{i+1},-b_{i+1}\}$ and $L_{i+1}=B(m;R_{i+1},S,T_{i+1})$.
We repeat the process until we obtain a connected graph $L_k=B(m;R_k,S,T_k)$ and set $L=L_k$. 

We now construct a hamilton cycle in the graph $L$.  First, we construct a hamilton cycle in each connected component  of $L_1$. To do that,  we use the construction from Lemma \ref{lem_BoPiZi2025_2ndpaper} if the components of $H$ are hamiltonian. 
Otherwise we use the constructions from Lemmas \ref{pro_method1_case}, \ref{pro_method1_case2} or \ref{pro_method1_case3}, as described in Section \ref{sec_construction}.
Note that we can apply the lemmas since $\lambda=\gcd(m,S,b_1)-1>0$, and, if $m$ is odd, also 
$\mu=\gcd(m,S)/\gcd(m,S,b_1)-1>0$
since $b_1$ is not divisible by $\gcd(m,S)$.
In fact, if $m$ is odd, then both parameters $\lambda$, $\mu$ are even and we apply Lemma \ref{pro_method1_case3}. 
For even $m$, we have also assumed that $m/\gcd(m, S)$ is odd.
In this case, if $b_1$ is odd, then $\lambda$ is even, $\mu$ is odd, and we apply Lemma \ref{pro_method1_case2}. 
If $b_1$ is even, then $\lambda$ is odd and we apply Lemma \ref{pro_method1_case2} or \ref{pro_method1_case} according to whether $b_1$ is divisible by the greatest power of $2$ dividing $\gcd(m,S)$ 
or not, respectively.

The constructions from Lemma \ref{lem_BoPiZi2025_2ndpaper} and Lemmas \ref{pro_method1_case}--\ref{pro_method1_case3} provide hamilton cycles containing at least two outer edges (see Remark \ref{rem_construction} and Figures \ref{fig_ham_cycle_lambda_mu_odd}--\ref{fig_ham_cycle_even}).
Therefore, we can apply Lemma  \ref{pro_construction_Bd1bis} to the components of $L_1$, and find a hamilton cycle in each component of $L_2$.
Since also the construction from Lemma \ref{pro_construction_Bd1bis} yields  hamilton cycles containing at least two outer edges, see Remark \ref{rem_construction}, we can apply Lemma  \ref{pro_construction_Bd1bis} to the components of $L_2$,  and find a hamilton cycle in each component of $L_3$. By repeatedly applying the construction from Lemma \ref{pro_construction_Bd1bis},  we find a hamilton cycle in $L_k=L$, and therefore in $G$,  since $L$ is a connected spanning subgraph of $G$. The assertion follows.
\end{proof}

\smallskip

We now establish Theorem \ref{th_a=b}.

\begin{proof}[Proof of Theorem \ref{th_a=b}]
Recall that $G=B(m; R, S, T)$ is a connected bicirculant graph with $|S|\ge 4$ and $\gcd(m, S)>1$. Moreover, the set $R\cap T$ contains an  element that is coprime to $\gcd(m, S)$. Let  $a$ be such an element. We construct a hamilton cycle in the connected spanning subgraph $G'=B(m; a, S, a)$ of $G$. This hamilton cycle is also a hamilton cycle in $G$.

We set $\lambda=\gcd(m, S)-1$ and denote with $H_i$, with $0\le i\le\lambda$, the connected components of $H(m; S)$; $H_0$ is the component containing $u_0$. The outer and inner vertices of the component $H_i$, with $0\le i\le\lambda$,  
are adjacent to the outer and inner vertices of $H_{i+1}$ through outer and inner edges of the same type $a$. 

By Lemma \ref{lemma_bicirculant_m<5}, we may assume that $m>5$. 
Therefore, by Lemma \ref{lemma_path_1},  graph $H_0$ has a hamilton path $P$  connecting an inner vertex $v_s$ to the vertex $u_0$,  with the property that at least one the neighbours of $v_s$,  say $u_z$, is such that the subpaths $v_s\,P\,u_z$, $u_z\,P\,u_0$ contain at least $6$ and $7$ vertices, respectively.  We denote the vertices of $P$ as in  Section \ref{sec_construction}. 

In analogy to the constructions in Lemma \ref{pro_method1_case}--\ref{pro_method1_case3}, we give a list of suitable paths
and then connect them to obtain a hamilton cycle in $G'$. 
We treat separately the cases in which $\lambda$ is even and odd.
\medskip

\noindent
\textbf{Case 1: $\lambda$ is even.}
The paths involved in the construction are the following:
the path $P$ in $H_0$,   the isolated  vertex $u^{\lambda-1}_0$, the edges  $v^{\lambda-1}_ku^{\lambda-1}_{z}$,  $v^{\lambda-1}_xu^{\lambda-1}_{x+1}$ with $0\le x\le k-1$, together the copy of the path $v_s\,P\,v_h$ in $H_{\lambda-1}$; we take
the edges $v^{\lambda}_xu^{\lambda}_{x}$ with $0\le x\le k$ and the copy of the path $u_zv_s\cup v_s\,P\,v_h$ in $H_{\lambda}$;
for $\lambda>2$ we also consider the isolated vertex $u^i_0$ together with a copy of the path $v_s\,P\,v_0$ in $H_i$, for every $1\le i\le\lambda-2$.

We connect the above paths and form a hamilton cycle in $G$ through the following edges:
the edges from the copy of $v_s$ in $H_i$ to the copy of $v_s$ in $H_{i+1}$ for every $0\le i\le\lambda-2$, $i$ even;
the edges from the copy of $v_0$ in $H_i$ to the copy of $v_0$ in $H_{i+1}$ for every $1\le i\le\lambda-1$, $i$ odd;
the edges from the copy of $u_0$ in $H_i$ to the copy of $u_0$ in $H_{i+1}$ for every $0\le i\le\lambda-1$;
the edges from the vertices $u^{\lambda-1}_x$, $v^{\lambda-1}_y$ in $H_{\lambda-1}$ to 
$u^{\lambda}_x$, $v^{\lambda}_y$, respectively, in $H_{\lambda}$, when  
$x=z$, $y=h$, and for $1\le x\le k$ and $1\le y\le k$.
\medskip

\noindent
\textbf{Case 2: $\lambda$ is odd.}
Take a copy of the path $P$ in $H_0$ and $H_{\lambda}$, a copy of the path $v_s\,P\,v_0$ together with a copy of the vertex $u_0$ in each $H_i$ with $1\le i\le\lambda-1$. Connect the paths  through the edges from $v_s$ in $H_{i}$ to $v_s$ in $H_{i+1}$ for $0\le i\le\lambda-1$, $i$ even, the edges from $v_0$ in $H_{i}$ to $v_0$ in $H_{i+1}$ for $1\le i\le\lambda-2$, $i$ odd, and the edges from $u_0$ in  $H_{i}$ to $u_0$ in $H_{i+1}$ for $0\le i \le\lambda-1$.
%
%
\end{proof}

\section{Concluding remarks}\label{sec_conclusion}

In this paper we have developed tools for the construction of a hamilton cycle in every bicirculant graph $B(m; R,  S, T)$ with $|S|\ge 4$, $\gcd(m,S)>1$  and such that $R\cup T$ contains at least one element that is not coprime to $\gcd(m, S)$ and is also not divisible by $\gcd(m, S)$ in the case where $m$ is odd. We construct a hamilton cycle in such a bicirculant  from hamilton paths in the connected components of its subgraph $H(m;S)$. This allows us to obtain results beyond those covered in our previous papers on this topic \cite{BoPiZi2025_2ndpaper,BoPiZi2025_3rdpaper}, where hamilton cycles in the subgraphs $H(m;S)$ were required. Additionally, we also show that in the case where all the elements of $R,T$ are coprime to $\gcd(m,S)$, but $R \cap T$ is nonempty, the graph $B(m;R,S,T)$ is hamiltonian. 
The accomplishments of this paper combined with the results in Section \ref{sec_preliminary} show that a bicirculant for which the existence of a hamilton cycle is not yet known satisfies one of the following conditions:

\begin{itemize}
\item $m$ is even,  $|S|\ge 4$,  $\gcd(m,S)>1$, $m/\gcd(m, S)$ is odd,  $m/2\not\in R,T$, $|R|=|T|>1$, $R\cap T=\emptyset$,  and every element in $R\cup T$ is coprime to $\gcd(m, S)$;

\item $m$ is odd,  $|S|\ge 4$ and either
\begin{itemize}
    \item $\gcd(m,S)=1$ and ($R \ne T$ or $R=T=\emptyset$) or 
    \item $\gcd(m,S)>1$, $R\cap T = \emptyset$ and every element in $R\cup T$ is either coprime to $\gcd(m, S)$ or divisible by $\gcd(m, S)$.
\end{itemize}
\end{itemize}
In particular, the existence of a hamilton cycle is not known for cyclic Haar graphs $H(m; S)$ with $m$ odd and $|S|\ge 4$.
%
By \cite[Theorem 1.5]{BoPiZi2025_3rdpaper},  we also know that every connected bicirculant graph with $|S|\ge 4$ is hamiltonian for even   $m< 9\, 240$,  and for odd $m< 3\,465$ or $m< 15\, 015$ according to whether the spokes induce a connected subgraph or not, respectively.  As remarked in \cite{BoPiZi2025_3rdpaper},  the first bicirculants for which the existence of a hamilton cycle is not known are,  for example, $B(9240; 2310, \{0, 70, 220, 154\}, 1155)$ if $m$ is even, 
$B(3465; 2310, \{0, 35, 77, 275\}, 1155)$ or $B(15015; 2310, \{0, 455, 3575, 1001\}, 1155)$ if $m$ is odd and
the spokes induce a connected subgraph or not, respectively. The results presented in this paper are not sufficient to establish the existence of a hamilton cycle in the above graphs,  since all elements of $R\cup T$ are coprime to $\gcd(m,S)$. It still remains a  challenge to find a hamilton cycle in these graphs because of their size.

We recall that the results in \cite{BoPiZi2025_3rdpaper} show that the existence of a hamilton cycle in every connected cyclic Haar graph implies that every connected bicirculant graph is hamiltonian. 
Proving that every connected cyclic Haar graph $H(m;S)$ with $m$ odd and $m \ge 3$ is hamiltonian would thus provide a complete positive solution to Conjecture \ref{conject}. 
However, we believe that this is a challenging problem and that further methods may need to be developed to obtain additional partial results and narrow down the family of open cases.

\section*{Acknowledgements}
Simona Bonvicini is a member of GNSAGA of Istituto Nazionale di Alta Matematica (INdAM).
The work of Toma\v{z} Pisanski is supported in part by the Slovenian Research and Innovation Agency
(research program P1-0284 and research project N1-0481).
Arjana \v Zitnik is supported in part by the Slovenian Research and Innovation  Agency (research program P1-0294 and research project J1-4351).
\bigskip

\noindent
\textbf{Data Availability Statement} No datasets were generated or analyzed during
the current study.

\section*{Declarations}

\noindent
\textbf{Conflict of interest} The authors declare that they have no conflict of interest.

\end{document}